\documentclass[12pt]{amsart}
\usepackage{amsmath, amsthm, amssymb}
\usepackage{fullpage}

\theoremstyle{plain}
\newtheorem{thm}{Theorem}

\newtheorem{lem}[thm]{Lemma}

\newtheorem*{MainLemma}{Main Lemma}
\newtheorem*{LovaszLemma}{Lov\'{a}sz's Decomposition Lemma}

\theoremstyle{definition}
\newtheorem{defn}{Definition}
\theoremstyle{remark}

\newcommand{\fancy}[1]{\mathcal{#1}}

\newcommand{\IN}{\mathbb{N}}

\newcommand{\G}{\fancy{G}}

\newcommand{\set}[1]{\left\{ #1 \right\}}

\newcommand{\card}[1]{\left|#1\right|}
\newcommand{\size}[1]{\left\Vert#1\right\Vert}

\newcommand{\func}[3]{#1\colon #2 \rightarrow #3}

\newcommand{\irange}[1]{\left[#1\right]}

\newcommand{\parens}[1]{\left( #1 \right)}

\newcommand{\DefinedAs}{\mathrel{\mathop:}=}

\title{A note on vertex partitions}
\author{Landon Rabern}

\begin{document}
\begin{abstract}
We prove a general lemma about partitioning the vertex set of a graph into subgraphs of bounded degree.  This lemma extends a sequence of results of Lov\'{a}sz, Catlin, Kostochka and Rabern.
\end{abstract}
\maketitle
\section{Introduction}
In the 1960's Lov\'{a}sz \cite{Lovasz} proved the following decomposition lemma for graphs by considering a partition minimizing a certain function.

\begin{LovaszLemma}
Let $G$ be a graph and $r_1, \ldots, r_k \in \IN$ such that $\sum_{i=1}^k r_i \geq \Delta(G) + 1 - k$. Then $V(G)$ can be partitioned into sets $V_1, \ldots, V_k$ such that $\Delta(G[V_i]) \leq r_i$ for each $i \in \irange{k}$.
\end{LovaszLemma}

A decade later, Catlin \cite{Catlin} showed that bumping the $\Delta(G) + 1$ to $\Delta(G) + 2$ allowed for shuffling vertices from one partition set to another and thereby proving stronger decomposition results. A few years later Kostochka \cite{Kostochka} modified Catlin's algorithm to show that every triangle-free graph $G$ can be colored with at most $\frac23 \Delta(G) + 2$ colors.  Around the same time, Mozhan \cite{Mozhan} used a different, but related, function minimization and vertex shuffling procedure to prove coloring results.  In \cite{rabern}, we generalized Kostochka's modification to prove the following.

\begin{lem}\label{LandonsLemma}
Let $G$ be a graph and $r_1, \ldots, r_k \in \IN$ such that $\sum_{i=1}^k r_i \geq \Delta(G) + 2 - k$. Then $V(G)$ can be partitioned into sets $V_1, \ldots, V_k$ such that $\Delta(G[V_i]) \leq r_i$ and $G[V_i]$ contains no non-complete $r_i$-regular components for each $i \in \irange{k}$.
\end{lem}

In fact, we proved a stronger lemma allowing us to forbid a larger class of components coming from any so-called \emph{$r$-permissible collection}.  The purpose of this note is to simplify and generalize this latter result.  The definition of an $r$-height function will be given in the following section.

\begin{MainLemma}
Let $G$ be a graph and $r_1, \ldots, r_k \in \IN$ such that $\sum_{i=1}^k r_i \geq \Delta(G) + 2 - k$. If $h_i$ is an $r_i$-height function for each $i \in \irange{k}$, then $V(G)$ can be partitioned into sets $V_1, \ldots, V_k$ such that for each $i \in \irange{k}$, $\Delta(G[V_i]) \leq r_i$ and $h_i(D) = 0$ for each component $D$ of $G[V_i]$.
\end{MainLemma}

\section{The proof}
Our notation follows Diestel \cite{Diestel} unless otherwise specified.  The natural numbers include zero; that is, $\IN \DefinedAs \set{0, 1, 2, 3, \ldots}$.  We also use the shorthand $\irange{k} \DefinedAs \set{1, 2, \ldots, k}$. Let $\G$ be the collection of all finite simple connected graphs.

\begin{defn}
For $\func{h}{\G}{\IN}$ and $G \in \G$, a vertex $x \in V(G)$ is called \emph{$h$-critical} in $G$ if $G - x \in \G$ and $h(G-x) < h(G)$.
\end{defn}

\begin{defn}
For $\func{h}{\G}{\IN}$ and $G \in \G$, a pair of vertices $\set{x,y} \subseteq V(G)$ is called an \emph{$h$-critical pair} in $G$ if $G - \set{x,y} \in \G$ and $x$ is $h$-critical in $G-y$ and $y$ is $h$-critical in $G-x$.
\end{defn}

\begin{defn}
For $r \in \IN$ a function $\func{h}{\G}{\IN}$ is called an \emph{$r$-height function} if it has each of the following properties:
\begin{enumerate}
\item if $h(G) > 0$, then $G$ contains an $h$-critical vertex $x$ with $d(x) \geq r$;
\item if $G \in \G$ and $x \in V(G)$ is $h$-critical with $d(x) \geq r$, then $h(G-x) = h(G) - 1$;
\item if $G \in \G$ and $x \in V(G)$ is $h$-critical with $d(x) \geq r$, then $G$ contains an $h$-critical vertex $y \not \in \set{x} \cup N(x)$ with $d(y) \geq r$;
\item if $G \in \G$ and $\set{x, y} \subseteq V(G)$ is an $h$-critical pair in $G$ with $d_{G-y}(x) \geq r$ and $d_{G-x}(y) \geq r$, then there exists $z \in N(x) \cap N(y)$ with $d(z) \geq r + 1$.
\end{enumerate}
\end{defn}

\noindent For $r \geq 2$, the function $\func{h}{\G}{\IN}$ which gives $1$ for all non-complete $r$-regular graphs and $0$ for everything else is an $r$-height function.  Applying the Main Lemma using this height function proves Lemma \ref{LandonsLemma}.

The proof of the Main Lemma uses ideas similar to those in \cite{Kostochka} and \cite{rabern}.  For a graph $G$, $x \in V(G)$ and $D \subseteq V(G)$ we use the notation $N_D(x) \DefinedAs N(x) \cap D$ and $d_D(x) \DefinedAs \card{N_D(x)}$. Let $\fancy{C}(G)$ be the components of $G$ and $c(G) \DefinedAs \card{\fancy{C}(G)}$.  If $\func{h}{\G}{\IN}$, we define $h$ for any graph as $h(G) \DefinedAs \sum_{D \in \fancy{C}(G)} h(D)$.

\begin{proof}[Proof of Main Lemma]
For a partition $P \DefinedAs \parens{V_1, \ldots, V_k}$ of $V(G)$ let

\[f(P) \DefinedAs \sum_{i=1}^k \parens{\size{G[V_i]} - r_i\card{V_i}},\]
\[c(P) \DefinedAs \sum_{i=1}^k c(G[V_i]),\]
\[h(P) \DefinedAs \sum_{i=1}^k h_i(G[V_i]).\]

\noindent Let $P \DefinedAs \parens{V_1, \ldots, V_k}$ be a partition of $V(G)$ minimizing $f(P)$, and subject to that $c(P)$, and subject to that $h(P)$.

Let $i \in \irange{k}$ and $x \in V_i$ with $d_{V_i}(x) \geq r_i$.  Since $\sum_{i=1}^k r_i \geq \Delta(G) + 2 - k$ there is some $j \neq i$ such that $d_{V_j}(x) \leq r_j$.  Moving $x$ from $V_i$ to $V_j$ gives a new partition $P^*$ with $f(P^*) \leq f(P)$.  Note that if $d_{V_i}(x) > r_i$ we would have $f(P^*) < f(P)$ contradicting the minimality of $P$. This proves that $\Delta(G[V_i]) \leq r_i$ for each $i \in \irange{k}$.

Now suppose that for some $i_1$ there is a component $A_1$ of $G[V_{i_1}]$ with $h_{i_1}(A_1) > 0$. Put $P_1 \DefinedAs P$ and $V_{1,i} \DefinedAs V_i$ for $i \in \irange{k}$. By property 1 of height functions, we have an $h_{i_1}$-critical vertex $x_1 \in V(A_1)$ with $d_{A_1}(x_1) \geq r_{i_1}$.  By the above we have $i_2 \neq i_1$ such that moving $x_1$ from $V_{1, i_1}$ to $V_{1, i_2}$ gives a new partition $P_2 \DefinedAs \parens{V_{2, 1}, V_{2,2}, \ldots, V_{2,k}}$ where $f(P_2) = f(P_1)$.  By the minimality of $c(P_1)$, $x_1$ is adjacent to only one component $C_2$ in $G[V_{2, i_2}]$. Let $A_2 \DefinedAs G[V(C_2) \cup \set{x_1}]$.  Since $x_1$ is $h_{i_1}$-critical, by the minimality of $h(P_1)$, it must be that $h_{i_2}(A_2) > h_{i_2}(C_2)$.  By property 2 of height functions we must have $h_{i_2}(A_2) = h_{i_2}(C_2) + 1$.  Hence $h(P_2)$ is still minimum.  Now, by property 3 of height functions, we have an $h_{i_2}$-critical vertex $x_2 \in V(A_2) - \parens{\set{x_1} \cup N_{A_2}(x_1)}$ with $d_{A_2}(x_2) \geq r_{i_2}$. 

Continue on this way to construct sequences $i_1, i_2, \ldots$, $A_1, A_2, \ldots$, $P_1, P_2, P_3, \ldots$ and $x_1, x_2, \ldots$.  Since $G$ is finite, at some point we will need to reuse a leftover component; that is, there is a smallest $t$ such that $A_{t + 1} - x_t = A_s - x_s$ for some $s < t$.  In particular, $\set{x_s, x_{t+1}}$ is an $h_{i_s}$-critical pair in  $Q \DefinedAs G\left[\set{x_{t+1}} \cup V(A_s)\right]$ where $d_{Q-x_{t+1}}(x_s) \geq r_{i_s}$ and $d_{Q-x_s}(x_{t+1}) \geq r_{i_s}$.  Thus, by property 4 of height functions, we have $z \in N_Q(x_s) \cap N_Q(x_{t+1})$ with $d_Q(z) \geq r_{i_s} + 1$.

We now modify $P_s$ to contradict the minimality of $f(P)$.  At step $t+1$,  $x_t$ was adjacent to exactly $r_{i_s}$ vertices in $V_{t+1, i_s}$. This is what allowed us to move $x_t$ into $V_{t+1, i_s}$.  Our goal is to modify $P_s$ so that we can move $x_t$ into the $i_s$ part without moving $x_s$ out. Since $z$ is adjacent to both $x_s$ and $x_t$, moving $z$ out of the $i_s$ part will then give us our desired contradiction.  

So, consider the set $X$ of vertices that could have been moved out of $V_{s, i_s}$ between step $s$ and step $t+1$; that is, $X \DefinedAs \set{x_{s+1}, x_{s+2}, \ldots, x_{t-1}} \cap V_{s, i_s}$.  For $x_j \in X$, since $d_{A_j}(x_j) \geq r_{i_s}$ and $x_j$ is not adjacent to $x_{j-1}$ we see that $d_{V_{s, i_s}}(x_j) \geq r_{i_s}$.  Similarly, $d_{V_{s, i_t}}(x_t) \geq r_{i_t}$. Also, by the minimality of $t$, $X$ is an independent set in $G$.  Thus we may move all elements of $X$ out of $V_{s, i_s}$ to get a new partition $P^* \DefinedAs \parens{V_{*, 1}, \ldots, V_{*, k}}$ with $f(P^*) = f(P)$. 

Since $x_t$ is adjacent to exactly $r_{i_s}$ vertices in $V_{t+1, i_s}$ and the only possible neighbors of $x_t$ that were moved out of $V_{s, i_s}$ between steps $s$ and $t+1$ are the elements of $X$, we see that $d_{V_{*, i_s}}(x_t) = r_{i_s}$.  Since $d_{V_{*, i_t}}(x_t) \geq r_{i_t}$ we can move $x_t$ from $V_{*, i_t}$ to $V_{*, i_s}$ to get a new partition $P^{**} \DefinedAs \parens{V_{**, 1}, \ldots, V_{**, k}}$ with $f(P^{**}) = f(P^*)$.  Now, recall that $z \in V_{**, i_s}$.  Since $z$ is adjacent to $x_t$ we have $d_{V_{**, i_s}}(z) \geq r_{i_s} + 1$.  Thus we may move $z$ out of $V_{**, i_s}$ to get a new partition $P^{***}$ with $f(P^{***}) < f(P^{**}) = f(P)$.  This contradicts the minimality of $f(P)$.
\end{proof}

\end{document}